\documentclass[a4paper,10pt]{article}
\usepackage{graphicx}
\begin{document}


\input amssym.tex
\parindent=0cm 
\hoffset=1.5truecm
\hsize=13truecm
\vsize=23truecm
\baselineskip 14pt


\def\Afigure{\cl{\includegraphics[width=.9\hsize]{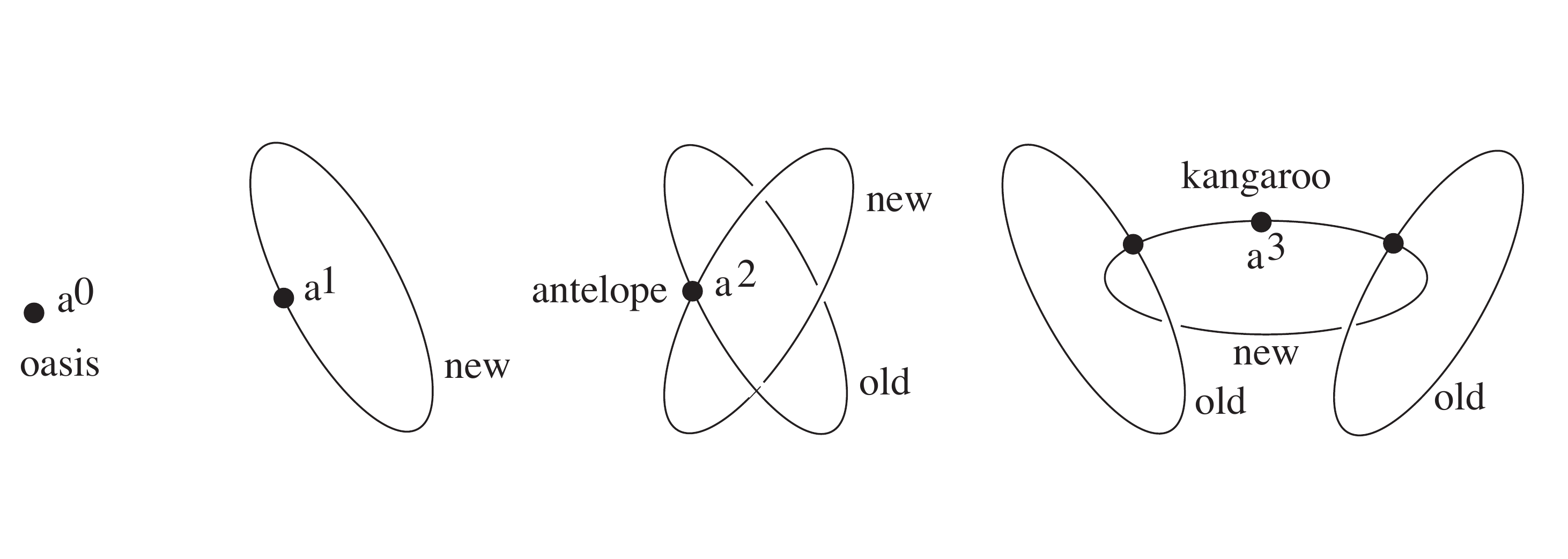}}}   


\font\Times=ptmr at 10pt
\font \Timesninepoint = ptmr at 9 pt
\font\ts=ptmr at 8 pt
\font\bf=ptmb at 10pt
\font\Bf=ptmb at 12pt
\font\bfit=ptmbi at 12pt
\font\ittenpoint=ptmri at 10pt
\Times

\def\litem{\par\noindent\hangindent=\parindent\ltextindent}
\def\ltextindent#1{\hbox to \hangindent{#1\hss}\ignorespaces}
\long\def\ignore#1\recognize{}
\long\def\<#1\>{\vbox{\leftskip=1truecm{\eightrm #1}}} 
\def\big{\bigskip}
\def\med{\medskip}
\def\meds{\medskip}
\def\hs{\hskip}
\def\vs{\vskip}
\def\hf{\hfill}
\def\no{\noindent}
\def\cl{\centerline}
\def\ds{\displaystyle}
\def\ol{\overline}
\def\sm{\setminus}
\def\ld{\ldots}
\def\cd{\cdot}
\def\to{,\ldots,}
\def\sub{\subseteq}
\def\{{\lbrace}
\def\}{\rbrace}
\def\iso{\cong}
\def\isom{\cong}
\def\congruent{\equiv}
\def\map{\rightarrow}
\def\inv{^{-1}}
\def\6{\partial}
\def\N{{\Bbb N}}
\def\C{{\Bbb C}}
\def\Z{{\Bbb Z}}
\def\Q{{\Bbb Q}}
\def\R{{\Bbb R}}
\def\A{{\Bbb A}}
\def\F{{\Bbb F}}
\def\H{{\Bbb H}}
\def\1{{\Bbb I}}
\def\O{{\cal O}}
\def\J{{\cal J}}
\def\weak{{\curlyvee}}
\def\ord{{\rm ord}}
\def\shade{{\rm shade}}
\def\min{{\rm min}}
\def\max{{\rm max}}
\def\a{\alpha}
\def\mm{m}
\def\abs#1{\vert#1\vert}
\def\nc{normal\ crossings\ }
\def\usc{upper\ semicontinuous\ }
\def\const{{\rm const}}
\def\r{r}
\def\q{q}
\def\t{t}
\def\olp#1{{\ol #1^p}}
\def\olc#1{{\ol #1^c}}
\def\minus {{\hbox {-}}}

\def\coeff{{\rm coeff}}
\def\jun{{\rm jun}}
\def\trans{{\rm trans}}
\def\contr{{\rm ctr}}
\def\label{{\rm lab}}
\def\owa{{\cal O}_{W,a}}
\def\red{{\rm red}}
\def\codim{{\rm codim}}
\def\edim{{\rm edim}}
\def\supp{\rm supp}
\def\o{{\bf o}}
\def\aa{{a}}
\def\strict{{st}}
\def\weak{{\curlyvee}}
\def\silent{{\tilde{}}}
\def\prior{{\circ}}
\def\-{{\hbox {-}}}
\def\={{\hbox {-:-}}}
\def\ord{{\rm ord}}
\def\min{{\rm min}}
\def\max{{\rm max}}
\def\a{\alpha}
\def\b{\beta}
\def\c{\gamma}
\def\abs#1{\vert#1\vert}
\def\nc{normal\ crossings\ }
\def\usc{upper\ semicontinuous\ }
\def\loc{{\rm loc}}
\def\ttop{{\rm top}}
\def\abs#1{\vert#1\vert}
\def\height{{\rm height}}

\def\vol{{\rm vol}}
\def\!{\bullet !\bullet}



\cl{\Bf  KANGAROO POINTS AND OBLIQUE POLYNOMIALS}\med

\cl{\Bf  IN RESOLUTION OF POSITIVE CHARACTERISTIC}\big\med

\cl{HERWIG HAUSER\footnote{\ts MSC-2000: 14B05, 14E15, 12D10.  Supported within the project P-18992 of the Austrian Science Fund FWF. The author thanks the Clay Institute for Mathematics at Cambridge and the Research Institute for Mathematical Science at Kyoto for their kind hospitality.}
}\big\big

The resolution of the singularities of algebraic varieties is still open in positive characteristic and arbitrary dimension. The inductive argument which works in characteristic zero fails for characteristic $p$. The main obstruction is the appearance of {\it kangaroo points}. These are points which show up in the course of a resolution process and for which the classical local resolution invariant {\it increases} (instead of decreasing). This destroys the induction.\med

Even though being very seldom, the phenomenon of kangaroo points had not been successfully overcome so as to allow a proof of resolution in characteristic $p$. So the author decided several years ago to investigate these points more closely. It seemed mandatory to analyze accurately what happened at these points. The results were astonishing: Aside of being very ``sporadic'', the singularities to be resolved showed a very specific pattern at kangaroo points. For instance, arithmetic conditions on the exceptional multiplicities had to be satisfied in order to allow an increase of the invariant. And, more strikingly, the coefficients of the involved polynomials seemed to play a decisive role -- a fact which contrasts our experience from characteristic zero. \med

The results of these investigations were written up around 2003 and assembled in the manuscript [Ha1], mostly for the author's personal reference. It was only circulated among the experts working in the field. After all, the characterization and classification of kangaroo points did not apparently show the way towards resolution in positive characteristic (even though a new proof for the surface case was found, as well as several other resolution strategies could be designed). Therefore, the paper [Ha1] was never published.\med

In fall 2008, Heisuke Hironaka gave several lectures at Harvard, the Clay Mathematics Institute at Cambridge and the Research Institute for Mathematical Sciences at Kyoto. There he presented a program about his view on how to attack the characterstic $p$ case. In the course of the lectures, Hironaka cited and used -- surprise -- the author's characterization of kangaroo points. He then claimed that this description {\it does indeed pave} the way towards resolution in positive characteristic. At the moment, no written confirmation of his claim is available, so its validity cannot be estimated yet. Nevertheless it seems worth to present the study of [Ha1] to a larger audience. \med

We will present here a concise overview on the theory of kangaroo points and their classification. The description is based on the notion of {\it oblique polynomials} treated in a later chapter. Proofs and more details can be found in the original paper [Ha1]. While adressing us mainly to algebraic geometers with some experience in resolution, we will add in footnotes explanations for readers which are curious but less familiar with the topic.
\med\med\goodbreak


{\bf Prelude for the non-expert reader.} Before getting into the actual material, let me briefly tell you what is resolution about and why it is important (and, also, why it is so fascinating). Readers acquainted with the subject may proceed directly to the next section. A system of polynomial equations in $n$ variables has a zeroset -- the associated {\it algebraic variety} $X$ -- whose structure can be quite complicated and mysterious. You may think of the real or complex solutions of an equation like $441(x^2y^2 + y^2z^2 + x^2z^2) = (1-x^2 - y^2 - z^2)^3$. The geometry of varieties shows all kind of local and global patterns which are difficult to guess from the equation. In particular, there will be {\it singularities}. These are the points where $X$ fails to be smooth (i.e., where $X$ is not a local manifold). At those points the Implicit Function Theorem IFT cannot be used to compute the nearby solutions. As a consequence, it is hard (also for computers) to describe correctly the local shape of the variety at its singular points.\med

Resolution of singularities is a method to understand where singularities come from, how they look like, and what is their internal structure. The idea is quite simple: When you take a submanifold $X$ of a high dimensional ambient space $M$ and then consider the image $X'$ of $X$ under the projection of the ambient space onto a smaller space $M'$, you most often create singularities on $X'$. The Klein bottle is smooth as a submanifold of $\R^4$, but there is no smooth realisation of it in $\R^3$. You necessarily have to accept self-intersections. Similarly, if you project a smooth space curve onto a plane in the direction of a tangent line at one of its points, the image curve will have singularities.\med

Which singular varieties can we obtain by such ``projections''? The answer is simple: {\it All!}\med

{\bf Theorem.} (Hironaka 1964) {\it Every algebraic variety over $\C$ is the image of a manifold under a suitable projection. A suitable manifold and map can be explicitly constructed (at least theoretically).}\med

For a geometer, this is quite amazing. For an algebraist, this is even more striking, since it means that it is possible to solve polynomial equations up to the Implicit Function Theorem. The applications of this result are numerous (it would be worth to list all theorems whose proofs rely on resolution). The reason is that, for smooth varieties, a lot of machinery is available to construct invariants and associated objects (zeta-functions, cohomology groups, characteristic classes, extensions of functions and differential forms, ...). As the projection map consists of a sequence of relatively simple maps (so called blowups), there is a good chance to carry these computations over to singular varieties. Which, in turn, is very helpful to understand them better.\med

Resolution is well established over fields of characteristic zero (with nowadays quite accessible proofs), but still unknown in positive characteristic (except for dimensions up to $3$). Why bothering about this? First, because (almost) everybody expects resolution to be true also in characteristic $p$. As the characteristic zero case was already a great piece of work (built on a truly beautiful concatanation of arguments), it is an intriguing challenge for the algebraic geometry community to find a proof which does not use the assumption of characteristic zero. But there is more to it: Many virtual results in number theory and arithmetic are just waiting to become true by having at the hand resolution in positive characteristic.\footnote{\ts De Jong's theory of {\it alterations}, valid  in arbitrary characteristic but slighty weaker than resolution, already produced a swarm of such results.} Again, it would be interesting to produce a list.\med

Another important feature of such a proof is our understanding of solving equations  in characteristic $p$. If we agree not to aim at one stroke solutions but to simplify the equation step by step (using for instance blowups) until we can see the solution (again, modulo IFT) there appears this delicate matter of understanding local coordinate changes in the presence of the Frobenius homomorphism. Phrased in very down to earth terms this means: {\it How do you measure whether a polynomial is, up to coordinate changes and up to adding $p$-th power polynomials, close or far from a monomial.} This may sound kind of silly, but be aware:  It is an extremely tough question (it resisted over 50 years) which lies at the very heart of the resolution of singularities in characteristic $p$. A meaningful proposal for such a measure (which should be compatible with blowups in a well defined sense) could break open the wall behind which we suspect to see a proof of resolution in positive characteristic. The rest would be mainly technicalities.\med

In the present article, we will see some of these ``elementary'' characteristic $p$ features, and we will make them very explicit. Of course it would be nice to have in parallel the conceptual counterparts of these constructions and phenomena, but this would require much more space and effort (for both the reader and the writer). As a consolation, the problems will be so concrete that everybody with a minimum talent in algebra will be tempted to attack them.\med

At the end of this paper, we briefly describe the actual state of research in resolution of singularities in positive characteristic (work of Hironaka, Villamayor, Kawanoue-Matsuki, W\l odarczyk, Cutkosky, Cossart-Piltant and others).

\med\med


{\bf Back to work.} Recall that the now classical resolution invariant in characteristic zero consists of a vector of integers whose components are orders of ideals in decreasing dimensions.\footnote{\ts For the basics on resolution, you may consult the survey [Ha2].} The ideals are the consecutive coefficient ideals in hypersurfaces of maximal contact, and the vector is considered with respect to the lexicographic ordering. It is then shown that this invariant drops under blowup in permissible centers.\footnote{\ts If you don't feel comfortable with blowups, we recommend [EiH] or [Ha5].} This allows to apply induction and to reduce by a sequence of blowups to the so called monomial case, for which an instant combinatorial description of the resolution is known. This program appears in different disguises in many places, see e.g. [Hi5, Vi1, BM, EV1, EH, W\l, Ko].\med

We start this note by reviewing the characteristic free version of the characteristic zero invariant of an ideal at a point as it was developed in [Ha1, EH]. For this definition, hypersurfaces of maximal contact (which need not exist in arbitrary characteristic) have to be replaced by hypersurfaces of {\it weak maximal contact}. These are defined as local regular hypersurfaces which {\it maximize} the order of the coefficient ideal of the given ideal (as hypersurfaces of maximal contact do), but whose transforms, in contrast,  are not required to contain along a sequence of blowups the points where the order of the original ideal remains constant. \med

Take then as resolution invariant the lexicographic vector consisting of the order of the ideal and of the orders of the iterated coefficient ideals with respect to such hypersurfaces. It turns out that the resulting vector (more precisely, its second component given by the order of the first coefficient ideal) may increase in positive characteristic under certain (permissible) blowups. The first examples of this phenomenon were observed by Abhyankar, Cossart, Moh and Seidenberg [Co, Mo, Se]. The increase destroys on first view any kind of induction. Moh succeeded to bound the maximal increase, but it was not yet possible to profit from this bound so as to save the induction argument (except for surfaces).  \med

We shall describe accurately the situations where an increase of the invariant occurs.  These are the {\it kangaroo points}.\footnote{\ts In [Hi1], kangaroo points run under the name of  metastatic points.} Actually, kangaroo points are very rare.  They are located at selected places of the exceptional divisor of a blowup, and lie above so called {\it antelope points}, which, in turn, can be completely classified.  \med

At an antelope point preceding a kangaroo point, three conditions hold: The residues modulo $p$ of the multiplicities of the exceptional components appearing in the defining equation must satisfy certain {\it arithmetic inequalities}, the order of the coefficient ideal of the equation must be {\it divisible} by the order of the equation, and {\it strong restrictions} on the (weighted) tangent cone of the defining equation are imposed (cf.÷  the theorem in section C). It turns out that  this tangent cone must be equal (up to multiplication by $p$-th powers) to an {\it oblique polynomial} in order to produce a kangaroo point after blowup. Oblique polynomials are characterized by a very particular behaviour under linear coordinate changes when considered up to addition of $p$-th powers. Fixing the exceptional multiplicities and the degree, both subject to the arithmetic and divisibility condition, it can be shown that there is precisely {\it one} oblique polynomial with these parameters (cf.÷ section E). \med


{\it Example.} This is the simplest example for the occurrence of a kangaroo point in a resolution process (cf.÷ section G for more details). Consider the following sequence of three point blowups in characteristic $2$, \med

\hskip 1.5cm $f^0 =x^2+ 1\cd(y^7+  yz^4)$ (oasis point $a^0$), \hf $(x,y,z)\map (xy,y,zy)$,\med

\hskip 1.5cm $f^1=x^2+y^3\cd (y^2+  z^4)$,  \hf $(x,y,z)\map (xz,yz,z)$,\med

\hskip 1.5cm $f^2=x^2+y^3z^3\cd(y^2+z^2)$ (antelope point $a^2$),  \hf $(x,y,z)\map (xz,yz+z,z)$,\med

\hskip 1.5cm $f^3=x^2+z^6\cd(y+1)^3((y+1)^2+1)$,\med

\hskip 1.95cm $ =x^2+z^6\cd(y^5+y^4+y^3+y^2)$  (kangaroo point $a^3$).\med

The oblique polynomial appears at the antelope point $a^2$ in the form $y^3z^3\cd(y^2+z^2)$. The kangaroo point is a uniquely specified point $a^3$ of the exceptional divisor of the third blowup. It lies off the transforms of the exceptional components produced by the first two blowups (see Figure 1). The coordinate change $x\map x+ yz^3$ at $a^3$ eliminates $y^2z^6$ and produces\med

\hskip 1.5cm $f^3 =x^2+z^6\cd(y^5+y^4+y^3)$.\med

The order of $f$ has remained constant equal to $2$ throughout. But the order of the coefficient ideal of $f$ (divided by the exceptional monomials) with respect to hypersurfaces of weak maximal contact has increased between $a^2$ and $a^3$. Namely, in $y^3z^3\cd(y^2+z^2)$ the monomial $y^3z^3$ is exceptional and the remaining factor $y^2+z^2$ has order $2$, whereas in $z^6\cd(y^5+y^4+y^3)$ the exceptional factor is $z^6$ and the remaining factor $y^5+y^4+y^3$ has order $3$.\med

\vs-.6cm
\Afigure  
\vs-.6cm
\cl{Figure 1: The configuration of kangaroo, antelope and oasis points.}
\med\med


For surfaces, it is possible to show that the resolution invariant {\it decreases in the long run}, i.e., that the occasional increases are compensated by decreases before and after them. A first method for proving this is developed in [Ha1] and will be sketched in section F below. A second, more general approach using the {\it bonus} of a singularity will appear in a forthcoming paper of Hauser, Wagner and Zeillinger [HWZ].\med \med

{\it Caveat.} It is human to try the resolution invariant from characteristic zero also in positive characteristic. After having observed that it may increase in special circumstances, it is also natural to study the cases where this actually happens. This will be done in this paper. The hope then is that the understanding of the obstruction may allow to overcome the increase either by extra arguments or by modifying the invariant thus yielding finally a complete resolution. This would be the conservative approach to characteristic  \med

But the geniune advance would consist in inventing a new invariant (which never increases). This would be -- in the simplest case -- a new measure which describes the ``distance'' of a polynomial to be a monomial (up to coordinate changes and multiplication by units in the formal power series ring). The classical recipe {\it factor from the polynomial the exceptional monomial and take the order of the remaining factor as invariant} seems to be just too crude in arbitrary characteristic.\med


{\it Acknowledgements.} The author is indebted to many people for sharing their ideas and insights with him, among them Heisuke Hironaka, Shreeram Abhyankar, Orlando Villamayor, Santiago Encinas, Ana Bravo, Gerd M\"uller, Josef Schicho, G\'abor Bodn\'ar, Edward Bierstone, Pierre Milman,  Dale Cutkosky, Jaroslav W\l odarczyk, Bernard Teissier, Vincent Cossart, Mark Spivakovsky, Hiraku Kawanoue, Kenji Matsuki, Li Li, Daniele Panazzolo, Anne Fr\"uhbis-Kr\"uger and J\'anos Koll\'ar. We thank Dominique Wagner and Eleonore Faber for a careful reading of the text and several substantial improvements, and Rocio Blanco for very helpful programming support. \med\med\goodbreak


{\bf A. The invariant.} We define only the first two components of the classical resolution invariant as these suffice for the phenomena to be described here. For an ideal sheaf $\J$ on a regular ambient scheme $W$ and a point $a\in W$ denote by $J=J_a$ the stalk of $\J$ at $a$. For convenience, we denote -- if appropriate -- by the same character $J$ the ideal generated in the completion $\widehat \O_{W,a}$ of the local ring $\O_{W,a}$.\footnote{\ts You may think here that $J$ is an ideal in a polynomial or formal power series ring.} For a local regular hypersurface $V$ in $W$ through $a$, the {\it coefficient ideal of $J$ in $V$} is defined as the ideal \med

\hs 3cm $\ds \coeff_VJ=\sum_{i=0}^{o-1}\ (a_{f,i}, f\in J)^{o!\over o-i}$,\med

where $o=\ord_aJ$ is the order of $J$ at $a$, $x=0$ is a local equation for $V$ and $f=\sum a_{f,i} x^i$ is the expansion of $f$ with respect to $x$, with coefficients $a_{f,i}\in \O_{V,a}$. Among the many variants of this definition in the literature, the given one suits best our purposes. More specifications appear in [EH].\med

In case $J$ is a principal ideal generated by one polynomial $f(x,y) =x^o+g(y)$ in $\A^{1+m}$ with variables $x$ and $y=(y_m\to y_1)$, the coefficient ideal of $J$ with respect to the hypersurface $x=0$ is simply the ideal in $\A^m$ generated by $g^{(o-1)!}$. The factorial is only needed to ensure integer exponents when $f$ has other $x$-terms.\med

The order of the coefficient ideal at $a$ depends on the choice of the hypersurface $V$, but remains unchanged under passing to the completions of the local rings. The supremum of these orders over all choices of local regular hypersurfaces $V$ through $a$ is a local invariant of $J$ at $a$ (i.e., by definition, only depends on the isomorphism class of the complete local ring $\widehat \O_{W,a}/J$).
This supremum is $\infty$ if and only if $J$ is {\it bold regular} at $a$, viz generated by a power of a parameter of $\widehat \O_{W,a}$ [EH]. If the supremum is $<\infty$ and hence a maximum, any hypersurface $V$ realizing this value is said to have {\it weak maximal contact with $J$ at $a$}. In characteristic zero, hypersurfaces of maximal contact have weak maximal contact [EH]. Moreover, their strict transforms under a permissible blowup $W'\map W$ contain all {\it equiconstant points} (= infinitely near points in $W'$), i.e., those points of the exceptional divisor where the order of the weak transform $J'$ of $J$ has remained constant (recall that this order cannot increase if $J$ has constant order along the center). \med

In arbitrary characteristic, the supremum of the orders of the coefficient ideal $\coeff_VJ$ for varying $V$ can be used to define the second component of the candidate resolution invariant of $J$ at $a$. If the supremum is $\infty$ and thus $J$ is bold regular, a resolution is already achieved locally at $a$, so we discard this case. We henceforth assume that the supremum is $<\infty$ and can thus be realized by the choice of a suitable hypersurface $V$. After factoring from the resulting coefficient ideal a suitable divisor one takes the order of the remaining factor as the second component of the invariant. More explicitly, let $D$ be a given normal crossings divisor in $W$ with defining ideal $I_WD$. We shall assume throughout that $\coeff_VJ$ factors for any chosen $V$ transversal to $D$ (in the sense of normal crossings) into a product of ideals \med

\hs 3cm $\coeff_VJ= I_V(D\cap V) \cd I_-$,\med

where $I_-$ is some ideal in $\widehat \O_{V,a}$ (this assumption is always realized in practice).  Then define the {\it shade of $J$ at $a$ with respect to $D$} as the maximum value $\shade_aJ$ of $\ord_aI_-$ over all choices of $V$ transversal to $D$.  In [Hi1], a similarly defined invariant is considered by Hironaka and called there the {\it residual order of $J$ at $a$}. As usual, questions of well-definedness and upper-semicontinuity have to be taken care of.\footnote{\ts \baselineskip 12pt  Semicontinuity works well if only closed points are considered. For arbitrary (i.e., non-closed) points, their appear pathologies which are described and studied by Hironaka [Hi1].}\med

Along a resolution process, $D$  will be supported by the exceptional locus at the respective stage.  It coincides with the second entry of  the {\it combinatorial handicap of a mobile} as defined in [EH]. At the beginning, or whenever $\ord_aJ$ has dropped, $D$ will be empty. If the order of $J$ has remained constant at a point $a'$ above $a$, the transform $D'$ of $D$ is defined as\med

\hs 3cm $D'=D^\weak + (\ord_a(D\cap V)+\shade_aJ-\ord_aJ)\cd Y'$,\med

where $Y'$ denotes the exceptional divisor of the last blowup, and $D^\weak$ the strict transform of $D$.\footnote{\ts We use here implicitly that $V$ and $Z$ are transversal to $D$. This is indeed the case in the resolution process of an ideal or scheme.} Note that $\ord_a(D\cap V)+\shade_aJ= \ord_a(\coeff_VJ)$. It follows from the transformation rule of $D$ that, under permissible blowup, the weak transform $J'$ of $J$ at an equiconstant point $a'$ above $a$ has as coefficient ideal $\coeff_{V'} J'$ in the strict transform $V'$ of $V$ an ideal which factors again into a product $I_{V'}(D'\cap V')\cd I_-'$, with $I_-'$ the weak transform $(I_-)^\weak$ of $I_-$. Here, it is assumed that the center $Z$ is contained in $V$. This is more delicate to achieve in positive characteristic, due to the example of Narasimhan where the singular locus of $J$ is not contained locally in any regular hypersurface [Na1, Na2, Mu]. It can, however, be realized by refining the usual stratification of the singular locus of $J$ through the local embedding dimension of this locus. \med

The {\it commutativity} of the passage to coefficient ideals with blowups can be subsumed as follows, cf.÷  [EH, Ha2]. Given a blowup with center $Z$ contained in the local hypersurface $V$ of $W$ locally at $a$ and transversal to $D$, we get for any equiconstant $a'$ in $W'$ above $a$ and $I_-'=(I_-)^\weak$ a commutative diagram\med

\hs 3cm $J'\hs .2cm \rightsquigarrow\hs .2cm \coeff_{V'}J'= I_{V'}(D'\cap V') \cd I_-'$\med

\hs 3.05cm $\downarrow$\hs 1.6 cm $\downarrow$  \med

\hs 3cm $J\hs .3cm \rightsquigarrow \hs .2cm \coeff_VJ\hs.2cm = I_V(D\cap V) \cd I_- $\med

 Here, the situation splits according to the characteristic: In characteristic zero, choosing for $V$ a hypersurface of maximal contact for $J$ at $a$, the strict transform $V'$ constitutes again a hypersurface of maximal contact for $J'$ at $a'$. In particular, both will have weak maximal contact so that the shades of $J$ and $J'$ are well-defined. In addition, $\shade_{a'}J'$ can be computed from $\shade_aJ$ by looking at the blowup $V'\map V$ with center $Z$ and the ideals $I_-$ and $I_-'$ (recall that $Z\subset V$ locally at $a$). As $\shade_aJ=\ord_aI_-$,  $\shade_{a'}J'=\ord_{a'}I_-'$ and $I_-'$ is the weak transform of $I_-$, it follows automatically that $\shade_{a'}J'\leq \shade_aJ$ (it is required here that the order of $I_-$ is constant along $Z$, a property that is achieved through the insertion of {\it companion ideals} as suggested by Villamayor, cf. [EV2, EH]). This makes the induction and the descent in dimension work. \med

In positive characteristic, it is in general not possible to choose a local hypersurface of maximal contact for $J$ at $a$. But a hypersurface of weak maximal contact will always exist, by definition. So choose one, say $V$. The good news is -- as already Zariski observed [Za] -- that the strict transform $V'$ of $V$ will contain all equiconstant points $a'$ of $J$ in the exceptional divisor $Y'$. The bad news is, as experimentation (or Moh's and Narasimhan's examples) show, that $V'$ need no longer have weak maximal contact with $J'$ at $a'$. Said differently, $V'$ need not maximize the order of the coefficient ideal of the weak transform $J'$ of $J$ at $a'$. One may have to choose a {\it new} hypersurface $U'$ at $a'$ to maximize this order. As Moh observed [Mo], there is still worse news, since the choice of $U'$ may produce a shade of $J'$ at $a'$ which is {\it strictly larger} than the shade of $J$ at $a$. This destroys the induction over the lexicographically ordered pair $(\ord_a(J), \shade_a(J))$. At least at first sight!\med

\med\goodbreak

{\bf B. Moh's bound.} In his paper on local uniformization, Moh investigates the possible increase of $\shade_aJ$ at equiconstant points $a'$ of $J$ in the purely inseparable case\med

\hs 3cm $f(x,y)= x^{p^e} + y^r\cd g(y)$,\med

with $\ord(y^r g)\geq p^e=\ord\ f$ and $e\geq 1$ (see [Mo]).\footnote{\ts Abhyankar has informed the author that he had communicated this observation to Moh.} Here, $V$ defined by $ x=x_n=0$ denotes a hypersurface of weak maximal contact for $f$ at $a=0$ in $W=\A^n$, $p$ is the characteristic of the (algebraically closed) ground field, and $y=(x_{n-1}\to x_1)$ denote further parameters so that $(x,y)$ form a complete parameter system of $R=\widehat \O_{\A^n,0}$. Moreover, $r\in\N^{n-1}$ is a multi-exponent whose entries are the multiplicities of the divisor $D\cap V$ at $0$, and $y_i=0$ defines an irreducible component of $D\cap V$ in $V$ for all $i$ for which $r_i>0$. All expressions take place in an \'etale neighborhood of $0$ in $\A^n$, so that $f$ and possible coordinate changes are considered as formal power series.  The shade of $f$ at $0$ with respect to the divisor $D$ defined by $y^r=0$ is given by $\shade_0 f=\ord_0g$, by the choice of $V$.\med


{\bf Proposition.} (Moh) {\it In the above situation, let $(W',a')\map (W,a)$ be a local blowup with smooth center $Z$ contained in the top locus of $f$ and transversal to $D$. Assume that $a'$ is an equiconstant point for $f$ at $a$, i.e., $\ord_{a'}f' =\ord_af=p^e$, where $f'$ denotes the weak (= strict) transform of $f$ at $a'$. Then \med

\hs 3cm $\shade_{a'}f' \leq \shade_af +p^{e-1}$.}\med


In case $e=1$, the inequality reads $\shade_{a'}f' \leq \shade_af +1$, which is not too bad, but still unpleasant. The short proof of Moh uses a nice trick with derivations, thus eliminating all $p$-th powers from $y^rg(y)$. He then briefly investigates the case where an increase of the shade indeed occurs, showing that in the next blowup the shade has to drop at least by $1$ (if $e=1$). This, as Moh underlines and we all know, does not suffice yet to make induction work.\med\med


{\bf C. Kangaroo points.} The purpose of the following paragraphs is to describe in compact form the classification of kangaroo points from [Ha1]. We recall: The shade of a polynomial $f$ at a point $a$ with respect to a normal crossings divisor $D$ is the {\it maximal} value of the order of its coefficient ideal minus the multiplicity of $D$ at $a$, the maximum taken over all choices of regular local hypersurfaces at $a$ transversal to $D$. A {\it kangaroo point} in a blowup $W'\map W$ with permissible center $Z$ and exceptional divisor $Y'$ is an equiconstant point $a'$ above  $a\in Z$ where the shade of $f$ with respect to $D$ has increased,\med

\hs 3cm  $\ord_{a'}f'=\ord_af$ \hs.3cm and \hs.3cm $\shade_{a'}f' > \shade_af$.\med

Here, $\shade_{a'}f'$ is taken with respect to the divisor \med

\hs 3cm $D'=D^\weak + (\ord_a(D\cap V)+\shade_af-\ord_af)\cd Y'$,\med

where $D^\weak$ denotes the strict transform of $D$ and $\ord_af=\ord_Zf$ holds by permissibility of $Z$. The point $a$ prior to a kangaroo point $a'$ is called {\it antelope point}. Note here that if $V$ and $V'$ are hypersurfaces of weak maximal contact  then $\ord_a (D\cap V) =\ord_aD$ and $\ord_{a'} (D'\cap V') =\ord_{a'}D'$ by transversality of $D$ and $D'$ with $V$ and $V'$.\med

For the ease of the exposition, we restrict to hypersurfaces in $\A^n=\A^{1+m}$ with purely inseparable equation $f(x,y)= x^p + y^r\cd g(y)$ of order $p$ at $0$ equal to the characteristic of the ground field, with exceptional multiplicities $r=(r_m\to r_1)\in\N^m$ and parameters $(x,y)=(x,y_m\to y_1)$. We shall work only at closed points and with formal power series. Moreover, we confine to point blowups, since these entail the most delicate problems. Most of the concepts and results go through for more general situations, cf. [Ha1]. For an integral vector $r\in \N^m$ and a number $c\in\N$, let $\phi_c(\r)$ denote the number of components of $\r$ which are not divisible by $c$, \med

\hs 3cm $\phi_c(\r)=\#\{i\leq m,\, r_i\not\congruent 0 \hbox { mod } c\}$.\med
 
Define $\olc\r=(\olc {r_m}\to \olc {r_1})$ as the vector
of the residues $0\leq \olc {\r_i}<c$ of the components of $r$ modulo $c$, and let $\abs r = r_m+\ldots +r_1$.\med

The next theorem characterizes kangaroo points.  This theorem is the result quoted in Hironaka's notes from September 2008 [Hi1, Prop.÷  13.1, Thm.÷  13.2]. We only formulate it here for purely inseparable polynomials of order $p$ at $0$. An appropriate extension also holds beyond the purely inseparable case and for blowups in positive dimensional centers, see  [Ha1, Thm.÷ 1, sec.÷ 5, and Thm.÷ 2, sec.÷ 12]. 
\med
\med


{\bf Kangaroo Theorem.} (Hauser) {\it Let $(W',a') \map (W,a)$ be a local point blowup of $W=\A^{1+m}$ with center $Z=\{a\}=\{0\}$ and exceptional divisor $Y'$. Let be given local coordinates $(x,y_m\to y_1)$ at $a$ so that $f(x,y)= x^p + y^r\cd g(y)\in \widehat\O_{W,a}$ has order $p$ and shade $\ord_ag$ with respect to the divisor $D$ defined by $y^r=0$. Let $f'$ be the strict transform of $f$ at $a'$. Set $D'=D^\weak+ (\abs r+\ord_ag - p)\cd Y'$ with $D^\weak$ the strict transform of $D$ in $W'$. Then, for $a'$ to be a kangaroo point for $f$, i.e., $\ord_{a'}f' =\ord_af$ and  $\shade_{a'}f'>\shade_af$, the following conditions must hold at $a$:}\med

(1) {\it The order $\abs r+\ord_ag$ of $y^rg(y)$ is a multiple of $p$.}\med

(2) {\it The exceptional multiplicities $r_i$ at $a$ satisfy }\med

\hs 3cm $\olp {r_m}+\ldots +\olp {r_1}\leq (\phi_p(\r)-1)\cd p$.\med

(3) {\it The point $a'$ is determined by the expansion of $f$ at $a$. It lies on none of the strict transforms of the exceptional components $y_i=0$ for which $r_i$ is not a multiple of $p$.}\med

(4) {\it The tangent cone of $g$ equals, up to linear coordinate changes and multiplication by $p$-th powers, a specific homogeneous polynomial, called {\rm oblique}, which is unique for each choice of $p$, $r$ and degree.}\footnote{\ts  The possibility of multiplication with $p$-th powers was not properly indicated in the original version of [Ha1] (though it was proven there).}\med


{\it Remarks.} (a) The necessity of condition (1) is easy to see and already appears in [Mo]. The arithmetic inequality in condition (2) is related to counting the number of $p$-multiples in convex polytopes and their $r$-translates in $\R^m$. It implies that {\it at least two} exponents $r_i$ must be prime to $p$. For surfaces ($m=2$), condition (2) reads $r_2, r_1 \not \congruent 0$ mod $p$ and $\ol r_2 + \ol r_1 \leq p$. Condition (3) implies that the reference point has to jump off from all exceptional components with $r_i\not \congruent 0 $ mod $ p$ in order to arrive at a kangaroo point. So it has to leave at least two exceptional components (cf.÷ Figure 1 from the introduction). This, together with the jump of $\shade_af$, justifies the naming of these points.  Condition (4) will be discussed in the example below and in section E on oblique polynomials. \med

(b) Conditions (1) to (4) are necessary for the occurrence of kangaroo points. They are also sufficient, up to the higher degree terms of $g$, in the following sense: In the transform $g'$ of $g$ the terms of $g$ of degree $>\ord_ag$ (i.e., not in the tangent cone) may have transformed into terms of degree smaller than the order of the transform $\ol g '$ of the tangent cone $\ol g$ of $g$. This signifies that $\ord_{a'}g'<\ord_{a'} \ol g'$. As $\ord_{a'}g'=\shade_{a'}f'$, $\ord_ag=\shade_af$ by definition, and $\ord_{a'} \ol g'\leq \shade_a\ol f +1=\shade_a f +1$ with $\ol f=x^p+y^r\ol g$ by Moh's bound applied to $\ol f$, the strict inequality $\shade_{a'}f'> \shade_af$ becomes impossible. The influence of the higher order terms of $g$ can be made quite explicit in concrete examples.\med

(c) We emphasize that the intricacy of the resolution in positive characteristic lies precisely in these higher order terms. Without them, $g$ is homogeneous (and thus equal to its tangent cone). In this case it is easy to make the order of $f$ drop below $p$ by suitable further blowups. But, in the general case, it seems to be tricky how to control $g$ beyond its tangent cone.\med


{\it Example.} For surfaces ($n=3$ and $m=2$), condition (2) reads $\ol r_2 + \ol r_1 \leq p$, provided that $r_2, r_1 >0$. In this case, there is an explicit description of the  tangent cone $P=\ol g$ of $g$ as indicated by condition (4): If ${k+r\choose k+1}$ is not a multiple of $p$ it has the form \med

\hs 1.6cm $P(y,z)=y^rz^s\cd\H_r^k(y,tz-y)=y^rz^s\cd\sum_{i=0}^k{k+r \choose i+r}y^{i}(tz-y)^{k-i}$\med

where $r=r_1$, $s=r_2$, $k=\ord_ag$ and $t$ is some non-zero constant in the ground field. The constant $t$ determines the location of $a'$ on the exceptional divisor $Y'$, and vice versa. The polynomials $\H_r^k(y,w)=\sum_{i=0}^k{k+r \choose i+r}
y^{i}w^{k-i}$ are called {\it hybrid polynomials of type $(r,k)$} in [Ha1]. Note that we can write $\H^k_r$ as \med

\hs 1.4cm $\H_r^k(y,w)= \sum_{i=0}^k{k+r \choose k-i}y^{i}w^{k-i}=$\med

\hs 2.8cm $= \sum_{i=0}^k {k+r \choose i}y^{k-i}w^{i}=$\med

\hs 2.8cm $=y^{-r}\cd \sum_{i=0}^k {k+r \choose i}y^{k+r-i}w^{i}=$\med


\hs 2.8cm $=\lfloor y^{-r}\cd \sum_{i=0}^{k+r} {k+r \choose i}y^{k+r-i}w^{i}\rfloor_{poly}=$\med

\hs 2.8cm $=\lfloor y^{-r}\cd (y+w)^{k+r}\rfloor_{poly}$,\med

where $\lfloor Q\rfloor_{poly}$ denotes those terms of the Laurent expansion of $Q$ which involve no monomials with negative exponents.  \med


Rocio Blanco observed that if ${k+r\choose k+1}$ is a multiple of $p$ the above polynomial $P=y^rz^s\cd\H_r^k(y,tz-y)$ is a $p$-th power and thus does not count. In this case one can use alternatively a description of the tangent cone of $g$ which is independent of the divisibility of ${k+r\choose k+1}$ by $p$ (cf. section E below):\med

\hs 2cm $P(y,z)= z^s\cd\int y^{r-1}(y-tz)^kdy=z^s\cd \sum_{i=0}^k \ (-1)^{k-i} {1\over r+i} \ y^{r+i}(tz)^{k-i}$,\med

the sum being taken over those $i$ for which $r+i$ is not divisible by $p$. Dominique Wagner showed that the two formulas for $P$ differ -- up to adding $p$-th powers -- by the scalar factor $(-1)^k{k+r\choose k+1}(k+1)$. This explains why the first formula requires that ${k+r\choose k+1}$ is prime to $p$.\med

Let us illustrate the dependence on $p$ in the case $p=k=2$, $r=s=2$, where the binomial coefficient ${k+r\choose k+1}={5\choose 3}=10$ is not prime to $p$. Indeed,\med

\hs 1.45cm $y^rz^s\cd\H_r^k(y,tz-y)=$\med

\hs 2.6cm $=y^3z^3\cd[{5\choose 3} (tz-y)^2 +{5\choose 4} y(tz-y)+{5\choose 5} y^2]=$\med

\hs 2.6cm $=y^3z^3\cd [10(tz-y)^2 +5 y(tz-y)+y^2]=$\med

\hs 2.6cm $=y^3z^3\cd [y(tz-y)+y^2]=$\med

\hs 2.6cm $=t y^4z^4$\med

is a $p$-th power (provided that $K$ is perfect) and thus does not count as oblique, whereas\med

\hs 1.45cm $P(y,z)= z^s\cd \int y^{r-1}(y-tz)^kdy=$\med

\hs 2.6cm $= z^3\cd \int y^2(y-tz)^2dy=$\med

\hs 2.6cm $= z^3\cd\int (y^4 + t^2y^2z^2) dy=$\med

\hs 2.6cm $= y^3z^3\cd (y^2+ t^2 z^2)$\med

produces an increase of the shade. In section E below, we characterize oblique polynomials in arbitrary dimension. \med

\med\med


{\bf D. Proof of the Kangaroo Theorem.} We indicate the main points of the argument for arbitrary polynomials $f$, i.e., in the case where $f$ is not necessarily purely inseparable. This makes things more complicated, but has the advantage to be generally applicable in a resolution process. The argument should be compared with the (much simpler) computation of oblique polynomials for the purely inseparable case which is given in section E. Along the way, one obtains an alternative proof of Moh's inequality. \med

It is convenient to work in the power series ring and to assume that $f$ is in Weierstrass form with respect to the variable $x$. It then suffices to consider a weighted homogeneous $f$  with respect to weights $(w,1\to 1)$ where $w\geq 1$ is the ratio between the order of $f$ and the order of its coefficient ideal with respect to $x$, say $w= \ord\ f /\ord\ \coeff_x(f) $. \med

First use the fact that the inequality $\olc {r_m}+\ldots +\olc {r_1}\leq (\phi_c(\r)-1)\cd c$ is equivalent to\med

\hs 3cm $ \lceil {\olc {r_m}\over c}\rceil +\ldots + \lceil {\olc {r_1}\over c}\rceil \leq
\lceil {\olc {r_m}+\ldots +\olc {r_1}\over c}\rceil$,\med

\no where $\lceil u\rceil$ denotes the smallest integer $\geq u$. This allows to count the lattice points which lie in certain integral simplices in $\R_+^n$  (called {\it zwickels} in [Ha1]) but do not belong to the sublattice $p\cd \Z^n$. The key step of the proof of the Kangaroo Theorem is then to establish the invertibility of the transformation matrix between the vectors of coefficients of polynomials with exponents in such zwickels under prescribed coordinate changes, the polynomials being always considered modulo $p$-th powers. For illustration, we reproduce the corresponding passage from section 11 of [Ha1]. \med\med


{\Timesninepoint \baselineskip 12pt 

Let $f(x,y)$ and $\tilde f(x,y)=f(x+\sum_\c h_\c y^\c,y+t y_m)$ be weighted homogeneous  polynomials  of weighted degree $e$ with respect to weights $(w,1\to1)$ on $(x,y)=(x,y_m\to y_1)$, where the sum $\sum_\c h_\c y^\c$ ranges over $\c\in\N^m$ with  ${\abs\c=w}$, and where $h_\c$ and the components of $t=(0, t_{m-1}\to t_1)$ belong to the ground field. Let $c=e/w$ be the order of $f$. Write  \med

\hs 3cm  $f(x,y)=\sum a_{k\a}x^ky^\a$ \ and \ $\tilde f(x,y)=\sum b_{l\b}(t)x^ly^\b$  \med

\no with $wk+\abs\a=wl+\abs\b=e$. We assume that $a_{c0}\neq 0$, i.e., that $x^c$ appears with non-zero coefficient, say $a_{c0}=1$. Let $V$ be the hypersurface in $W=\A^n$ defined by $x=0$.  
Let \med

\hs 1.5cm $L_c=\{(k,\a)\in \N^{1+m},\, k<c\}\map \Q^m:\, (k,\a)\map {c\over c-k}\cd\a$
\med

\no be the map projecting elements $(k,\a)$ of the layer
$L_c$ in $\N^{1+m}$ to elements of $\Q^m$. The center of the projection is the
point $(c,0\to 0)$.\med

Let  $\q\in \N^m$ with $\abs\q=q_m+\ldots + q_1\leq e$ be fixed. Define the {\it upper zwickel} $Z(\q)$ in $\N^{1+m}$ as the set of points $(k,\a)$  with   $0\leq k\leq c$, $wk+\abs\a=e$ and projection ${c\over c-k}\cd\a \geq_{cp}\q$, denoting by $\geq_{cp}$ the componentwise order.  Thus $Z(\q)$ is given by \med

\hs 2cm $Z(\q): \, wk+\abs\a=e$ \ and \ $\a \geq_{cp}\lceil{c-k\over
c}\cd(\q_m\to \q_1)\rceil$.\med

\no Let us fix a decomposition $\q=\r+\ell \in \N^m$ with $\r=(\q_m\to \q_{j+1},0\to 0)$ and $\ell=(0\to 0,\q_j\to \q_1)$ for some index $j$ between $m-1$ and $0$. Define the {\it lower zwickel} $Y(\r,\ell)$ in $\N^{1+m}$ as the set of points $(k,\b)$ in
$\N^{1+m}$ with  $0\leq k\leq c$, $wk+\abs\b=e$ and projection ${c\over c-k}\cd\b
\geq_{cp}(\abs{\r},0\to 0, \ell)$. Thus $Y(\r,\ell)$ is given by  \med

\hs .5cm  $Y(\r,\ell): \, wk+\abs\b=e$ \ and \ $\b\geq_{cp} \lceil({c-k\over
c}\cd\abs{\r},0\to 0, {c-k\over c}\cd\q_j\to {c-k\over c}\cd\q_1)\rceil$.\med

\no For $j=m-1$ and hence $\r=(q_m,0\to 0)$ and  $\ell=(0,q_{m-1}\to q_1)$ we have
$Z(\q)=Y(\r,\ell)$. In general, the two zwickels are different.\med

For any $\r$ and $\ell$ and $0\leq k\leq e/w=c$ the slice
\med

\hs 2cm $Y(\r,\ell)(k)= \{(k,\b)\in Y(\r,\ell)\}=Y(r,\ell)\cap (\{k\}\times\N^m)$ \med

\no has at least as many elements as the slice \med

\hs 2cm $Z(\q)(k)=\{(k,\a)\in Z(q)\}=Z(q)\cap (\{k\}\times\N^m)$. \med

\no This holds for $k=0$, by definition of $Z(q)$ and $Y(r,\ell)$. For
arbitrary $k$, the inequality $\lceil{c-k\over c}\cd\abs{\r}\rceil  \leq \abs{\lceil {c-k\over
c}\cd\r\rceil}$ implies that the condition  \med

\hs 1cm  $wk+\abs\b=e$ \ and \ $\b\geq_{cp} (\abs{\lceil {c-k\over
c}\cd \r\rceil },0\to 0, \lceil{c-k\over c}\cd \q_j\rceil\to \lceil{c-k\over c}\cd
\q_1\rceil)$\med

\no is more restrictive than the condition \med

\hs 1cm  $wk+\abs\b=e$ \ and \ $\b\geq_{cp} (\lceil {c-k\over
c}\cd \abs{\r}\rceil ,0\to 0, \lceil{c-k\over c}\cd \q_j\rceil\to\lceil{c-k\over c}\cd
\q_1\rceil)$\med

\no defining $Y(\r,\ell)(k)$. For each $k$, the set of pairs $k,\b$ satisfying the first
condition has as many elements as $Z(q)(k)$ because $\abs r + q_j +\ldots +q_1=\abs q$. The
claim follows.\med


It is immediate that $y^\q$ is a factor of $\coeff_Vf$ if and only if $f - x^c$ has all exponents in the upper zwickel $Z=Z(\q)$ with $\q\in\N^m$, and $\coeff_V\tilde f$ has order $> e-\abs r$ in $z=(y_{m-1}\to y_1)$  if and only if all
coefficients of the monomials of $\tilde f- x^c$ with exponent in the lower zwickel $Y(\r,\ell)$ are zero. \med

Write elements $\b\in \N^m$ as $(\b_m,\b^\minus )$ where
$\b^\minus =(\b_{m-1}\to \b_1)\in\N^{m-1}$. Let $Y^*(\r,\ell)$ be the subset of $Y(\r,\ell)$ of elements $(k,\b)\in \N^{1+m}$ given by\med

\hs 3cm $\abs{\b^\minus } \leq e-wk-\lceil{c-k\over c}\cd \abs r\rceil$,\med

\hs 3cm $\b^\minus   \geq_{cp}  \lceil{c-k\over c}\cd(0\to 0, \q_j\to q_1)\rceil$.\med

\no By definition, for each $k$, the slice $Y^*(\r,\ell)(k)$ has the same cardinality as
the slice $Z(\q)(k)$ of the upper zwickel $Z(\q)$. For $\a$ and $\delta$ in $\Z^m$ set
${\a\choose\delta}=\prod_i {\a_i\choose\delta_i}$ where ${\a_i\choose\delta_i}$ is zero if
$\a_i<\delta_i$ or $\delta_i<0$. For $\Gamma$ a subset of $\N^m$, define for $k\in \N$ and
$\lambda=(\lambda_\c)_{\c\in \Gamma}\in \N^\Gamma$ the alternate binomial coefficient
\medskip

\hs 3cm $ [{k\choose \lambda}]=\prod_{\c\in\Gamma} {k - \abs {\lambda}^\c \choose
\lambda_\c}$ \ with \ $\abs{\lambda}^\c=  \sum_{\varepsilon\in\Gamma, \varepsilon<_{lex}\c}
\lambda_\varepsilon $.\medskip

\no Let $\Gamma\subset\N^m$ be the set of $\c\in\N^m$ with $\abs{\c}=w$ and write
$h=(h_\c)_{\c\in \Gamma}$.  Set
${\lambda\cd\Gamma}=\sum_{\c\in\Gamma}\lambda_\c \cd\c \in\N^m$ and fix
$t=(0,t_{m-1}\to t_{j+1},0\to 0)$.  We then have [Ha1, Prop.÷ 1, sec.÷ 11]: 
\med\med\goodbreak

{\bf Proposition.÷} {\ittenpoint Let $f(x,y)=\sum a_{k\a}x^ky^\a$ and $\tilde f(x,y)=
f(x+\sum_{\c\in\Gamma} h_\c y^\c,y+t y_m)=\sum b_{l\b}(t)x^ly^\b$ be weighted homogeneous polynomials with respect to weights $(w,1\to 1)$ as above. Fix $q=r+\ell\in\N^m$ with zwickels $Z(\q)$ and  $Y^*(\r,\ell)\subset Y(r,\ell)$.}\medskip

\no (1) {\ittenpoint The transformation matrix $A=(A_{k\a,l\b})$ from the coefficients
$a_{k\a}$ of $f$ to the coefficients $b_{\l\b}(t)$ of $\tilde f$ is given by \med

\hs 2cm $A_{k\a,l\b}=\sum\limits_{\lambda\in\N^\Gamma, \abs \lambda =k-l} {k\choose l}[{k-l\choose \lambda}]{\a\choose \delta_{\a\b \lambda}}\cd  h^\lambda\cd  \t^{\a-\delta_{\a\b \lambda}},$\med

\no where $\delta_{\a\b \lambda}=(\a_m,\b^\minus  -({\lambda\cd\Gamma})^\minus )\in \N^m$ and $h^\c=\Pi_\c h_\c^{\lambda_\c}$.}\med


\no (2) {\ittenpoint The quadratic submatrix $A^\square=(A_{k\a,l\b})$ of $A$ with $(k\a,l\b)$ ranging in $Z(\q)\times Y^*(\r,\ell)$ has determinant \ $t^{\rho(Z,Y^*(\r,\ell))}$ where $\rho(Z,Y^*(\r,\ell))$ is a vector in $\N^{m-1}$ independent of $h=(h_\c)_{\c\in\Gamma}$ with $\rho_m=0$ and $\rho_j=\cdots=\rho_1=0$.}\med 


\no (3) {\ittenpoint Assume that $f$ has support in $Z(\q)$. If $t_{m-1}\to t_{j+1}$ are non zero, the coefficients $b_{l\b}$ of $\tilde f$ in the lower zwickel $Y(\r,\ell)$ determine all coefficients of $f$.  }\med

}

This ends the excerpt from [Ha1] about the proof of the Kangaroo Theorem. Actually, the assertions of the theorem are rather straightforward consequences of the above proposition: Inverting the transformation matrix between the coefficients vectors of the polynomials allows to determine the tangent cone of $g$ as alluded to in assertion (4) of the theorem. As for the proof of the proposition itself, the formula from (1) is an exercise in binomial expansion, assertion (2) is tricky and relies on a special numbering of the lattice points in zwickels in order to make the matrix block-diagonal, and (3) follows rather quickly from (2). \med
\med\goodbreak


{\bf E. Oblique polynomials.} We now describe the tangent cone of the polynomials $g$ appearing in $f=x^p+y^rg(y)$ at antelope points preceding a kangaroo point. In [Ha1], the uniqueness assertion (4) in the theorem above was established for the tangent cone of arbitrary hypersurfaces of order $p$, and oblique polynomials were characterized in various specific situations. In [Hi1], a general description of oblique polynomials is given, and Schicho found independently a similar formula. Below we combine all viewpoints to a conjoint presentation. \med

Fix variables $y=(y_m\to y_1)$. Set $\ell=m-1$, and let $p$ be the characteristic of the ground field $K$. A non-zero polynomial $P=y^r g(y)$ with $r\in \N^m$ and $g$ homogeneous of degree $k$ is called {\it oblique with parameters} $p$, $r$ and $k$ if $P$ has no non-trivial $p$-th power polynomial factor and if there is a vector $t=(0,t_\ell\to t_1)\in (K^*)^m$ so that the polynomial $P^+(y)=(y+ty_m)^r  g(y+ty_m)$ has, after deleting all $p$-th power monomials from it, order $k+1$ with respect to the variables $y_\ell\to y_1$. Without loss of generality, the vector $t$ can and will be taken equal to $(0,1\to 1)$. We shall write $\ord^p_{z}P^+$ to denote the order of $P^+$ with respect to $z=(y_\ell\to y_1)$ modulo $p$-th powers. \med


{\it Example.} Take $m=2$, $p=2$ and $P(y)=y_2y_1(y_2^2+y_1^2)$ with $k=2$. Then $P^+(y)=P(y_2,y_1+y_2)= y_2y_1^2(y_1+y_2)$ has modulo squares order $3$ with respect to $y_1$.\med


It is checked by computation that the condition $\ord^p_{z}P^+ \geq k+1$ on $P^+$ is a prerequisite for the occurence of a kangaroo point as in the theorem. The result of Moh implies $\ord^p_z P^+\leq k+1$, so that equality must hold. Condition (4) of the theorem tells us that there is, up to addition of $p$-th powers, {\it at most one} oblique polynomial for each choice of the parameters $p$, $r$ and $k$. In order that $P$ is indeed oblique it is then also necessary that the degree of $P$ is a multiple of $p$ and that $r$ satisfies $\olp {r_m}+\ldots +\olp {r_1}\leq (\phi_p(\r)-1)\cd p$ (again by the theorem). \med

The following trick for characterizing oblique polynomials appears in [Ha1] for surfaces and is extended in [Hi1] to arbitrary dimension.  We dehomogenize $P$ with respect to $y_m$. This clearly preserves $p$-th powers. Moreover, when applied to monomials of total degree divisible by $p$ (as is the case for the monomials of the expansion of $P$), the dehomogenization creates no new $p$-th powers. It is thus an ``authentic'' transformation in our context, i.e., the characterization of oblique polynomials can be transcribed entirely to the dehomogenized situation.  Setting $y_m=1$ and $z=(y_\ell\to y_1)$ we get $Q(z)=P(1,z)=z^s\cd h(z)$ with $s=(r_\ell\to r_1)\in\N^\ell$ and $h(z)=g(1,z)$ a polynomial of degree $\leq k$. The translated polynomial is $Q^+(z)= Q(z+\1)=(z+\1)^s\cd h(z+\1)$, where $\1=(1\to 1)\in \N^\ell$.  The condition $\ord^p_z P^+\geq k+1$ now reads $\ord^p_z Q^+\geq k+1$ or, equivalently, $Q^+\in \langle z_\ell\to z_1\rangle^{k+1}+  K[z^p]$. Let us write this as\med

\hs 3cm $(z+\1)^s\cd h(z+\1) - v(z)^p\in \langle z_\ell\to z_1\rangle^{k+1}$\med

for some polynomial $v\in K[z]$. As $h$ has degree $\leq k$, the polynomial $v$ cannot be zero. In addition, we see that the condition $\ord^p_z Q^+\geq k+1$ is stable under multiplication with homogeneous $p$-th power polynomials $w(z)$, in the sense that $\ord^p_z\ (w^p \cd Q^+)\geq k+1+p\cd \deg w$.
Using that $(z+\1)^s$ is invertible in the completion $K[[z]]$ we get\med

\hs 3cm $ h(z+\1) =\lfloor (z+\1)^{-s}\cd v(z)^p\rfloor_k$, \med

where $\lfloor u(z) \rfloor_k$ denotes the $k$-jet (= expansion up to degree $k$) of a formal power series $u(z)$. From Moh's inequality we know that $(z+\1)^s\cd h(z+\1) - v(z)^p$ cannot belong to $\langle z_\ell\to z_1\rangle^{k+2}$. Therefore, in case that $v(z)$ is a constant, the homogeneous form of degree $k+1$ in $(z+\1)^{-s}$ must be non-zero. This form equals $\sum_{\alpha \in \N^\ell, \abs\alpha =k+1} {-s\choose \alpha} z^\alpha$. We conclude that if all ${-s\choose \alpha}$ with $\abs\alpha =k+1$ are zero in $K$, then $v$ was not a constant.\footnote{\ts The converse need not hold, see the example.} Inverting the translation $\tau(z)=z+\1$ we get the following formula for the dehomogenized tangent cone at antelope points preceding kangaroo points,\med

\hs 3cm $ z^s\cd h(z) =z^s \cd \tau^{\-1}{\{\lfloor (z+\1)^{-s}\cd v(z)^p\rfloor_k \}}$. \med

The homogenization of this polynomial with respect to $y_m$ followed by the multiplication with $y_m^{r_m}$ then yields the actual oblique polynomial $P(y)=y^r  g(y)$.\med

{\it Example.} In the example $P(y)=y_2^3y_1^3(y_2^2+y_1^2)$ from the beginning we have characteristic $p=2$, exponents $r_2=r_1=3$ and degree $k=2$. Therefore $\ell =1$ and $s=3$, which yields a binomial coefficient ${-3\choose \alpha}= {-3\choose 3} = -10$ equal to $0$ in $K$. Indeed, $P$ has as non-monomial factor $g(y)$ the square $(y_2+y_1)^2$. In the example $P(y)=y_2y_1(y_2^2+y_1^2)$ from above with $r_2=r_1=s=1$, the polynomial $g$ is again a square, even though ${-s\choose \alpha}={-1\choose 3}=-1$ is non-zero in $K$.

\med\med\goodbreak

{\bf F. Resolution of surfaces.}  In the surface case, there are several ways to overcome (or avoid) the obstruction produced by the appearance of kangaroo points. The first proof of surface resolution in positive characteristic is due to Abhyankar, using commutative algebra and field theory [Ab1]. Resolution invariants for surfaces then appear, at least implicitly, in his later work on resolution of three-folds. In [Hi4], Hironaka proposes an explicit invariant for the embedded resolution of surfaces in three-space (see [Ha3] for its concise definition). It is not clear how to extend this invariant to higher dimensions. \med

In [Ha1], it is shown for surfaces that during the blowups prior to the jump at a kangaroo point the shade must have decreased at least by $2$ (with one minor exception) and thus makes up for the later increase at the kangaroo point. To be more precise, given a sequence of point blowups in a three dimensional ambient space for which the subsequent centers are equiconstant points for some $f$, call {\it antelope point} the point $a$ immediately prior to a kangaroo point $a'$, and {\it oasis point} the last point $a^\circ$ below $a$ where none of the exceptional components through $a$ has appeared yet. The following is then a nice exercise:\med

{\bf Fact.} {\it The shade of $f$ drops between the oasis point $a^\circ$ and the antelope point $a$ of a kangaroo point $a'$ at least to the integer part of its half, \med

\hs 3cm $ \shade_{a}f \leq \lfloor {1\over 2}\cd  \shade_{a^\circ}f^\circ\rfloor$.}\med\goodbreak

In the purely inseparable case of an equation of order equal to the characteristic, this decrease thus dominates the later increase of the shade by $1$ except for the case $\shade_{a^\circ}f^\circ=2$ which is easy to handle separately and will be left to the reader. It seems challenging to establish a similar statement for singular three-folds in four-space.  \med

In [HWZ], we proceed somewhat differently by considering also blowups after the occurence of a kangaroo point. A detailed analysis shows that when taking three blowups together (the one between the antelope and the kangaroo point, and two more afterwards), the shade always either decreases in total, or, if it remains constant, an auxiliary secondary shade drops. This shade can again be interpreted as the order of a suitable coefficient ideal (now in just one variable), made coordinate independent by maximizing it over all choices of hypersurfaces inside the chosen hypersurface of weak maximal contact.  \med

The cute thing is that one can subtract, following an idea of Dominik Zeillinger [Ze] which was made precise and worked out by Wagner, a correction term from the shade which eliminates the increases without creating new increases at other blowups. This correction term, called the {\it bonus}, is defined in a subtle way according to the internal structure of the defining equation. It is mostly zero, takes at kangaroo points a value between $1$ and $2$, and in certain well defined situations a value between $1/2$ and $1$.\med

This bonus allows to define an invariant -- a triple consisting of the order, the modified shade and the secondary shade -- which now drops lexicographically after {\it each} blowup. The bonus is defined with respect to a {\it local flag} as defined in [Ha4]. Flags break symmetries and are stable under blowup (in a precise sense) and thus allow to define the bonus at any stage of the resolution process. We refer to [HWZ] for the details, as well as for the definition of an alternative invariant, the {\it height}, which is even simpler to use for the required induction. It profits much more from the flag than the shade and allows a simpler definition of the bonus. The invariant built from the height yields a quite systematic induction argument which may serve as a testing ground for the embedded resolution of singular three-folds.\med\med


{\bf G. Example of kangaroo point.} We comment on the example from the introduction on the occurence of kangaroo points. We are in characteristic $2$ and consider three point blowups, \med

\hskip 1.5cm $f^0 =x^2+ 1\cd(y^7+  yz^4)$ (oasis point $a^0$),  \hf $(x,y,z)\map (xy,y,zy)$,\med

\hskip 1.5cm $f^1=x^2+y^3\cd (y^2+  z^4)$, \hf $(x,y,z)\map (xz,yz,z)$,\med

\hskip 1.5cm $f^2=x^2+y^3z^3\cd(y^2+z^2)$ (antelope point $a^2$), \hf $(x,y,z)\map (xz,yz+z,z)$,\med

\hskip 1.5cm $f^3=x^2+z^6\cd(y+1)^3((y+1)^2+1)$, \med

\hskip 1.95cm $ =x^2+z^6\cd(y^5+y^4+y^3+y^2)$  (kangaroo point $a^3$).\med

The oblique polynomial appears at $a^2$ in the form $y^3z^3\cd(y^2+z^2)$. The kangaroo point $a^3$ occurs in the third blowup and is the unique equiconstant point of the exceptional divisor where the shade of $f$ increases. It lies off the transforms of the exceptional components produced by the first two blowups. In $W^3$, the strict transform of the hypersurface $x=0$ in $W^2$ has no longer weak maximal contact. The coordinate change $x\map x+ yz^3$ is needed to realize the shade, yielding in the new coordinates the expansion $f^3 =x^2+z^6\cd(y^5+y^4+y^3)$ with $\shade_{a^3} f^3=3> 2=\shade_{a^2}f^2$. Observe that the shade drops between the oasis and antelope point by $3$.

\med\med\goodbreak


{\bf H. Bibliographic comments.} We briefly relate the contents of this note to the existing literature on resolution in positive characteristic. The arithmetic condition $\olp {r_m}+\ldots +\olp {r_1}\leq (\phi_p(\r)-1)\cd p$ on the exceptional multiplicities at an antelope point appears in a different perspective also in the work of Abhyankar on good points [Ab2]. \med

There are important recent results and proofs of Cutkosky and Cossart-Piltant for the non-embedded resolution of three-folds in positive characteristic [Cu, CP]. Cutkosky reduces Abhyankar's proof  (over 500 pages) of resolution in characteristic $> 5$ to some forty pages, Cossart and Piltant establish the result with considerably more effort for arbitary fields. Both proofs use substantially the embedded resolution of surfaces (built on the invariant from [Hi4]), but they do not provide {\it embedded} resolution of three-folds. \med

As for dimension $n$, Hironaka develops in [Hi1, Hi2, Hi3] an elaborate machinery of differential operators in arbitrary characteristic in order to construct generalizations of hypersurfaces of maximal contact by allowing primitive elements as defining equations. The main difficulty is thus reduced to the purely inseparable case and  metastatic points, which precisely correspond to our kangaroo points. Hironaka then asserts that this type of singularities can be resolved directly.  There is no written proof of this available yet.\med 

There is a novel and impressive approach to resolution by Villamayor and his collaborators Benito, Bravo and Encinas [Vi2, Vi3, BV, EV3]. It is based on projections instead of restrictions for the descent in dimension.  A substitute for coefficient ideals is constructed via Rees algebras and differential operators, called elimination algebras. It provides a new resolution invariant for characteristic $p$ (which coincides with the classical one in zero characteristic). All the necessary properties are proven. This allows to reduce by blowups to a so called {\it monomial case} (which, however, seems to be still unsolved, and could be much  more intricate than the classical monomial case). \med

In a somewhat different vein, Kawanoue and Matsuki have announced a very promising program for resolution in arbitrary characteristic and dimension [Ka, MK]. Again, they use differential operators to define a suitable resolution invariant and then show its upper semicontinuity. The termination of the resulting algorithm seems not to be ensured yet.\med

W\l odarczyk has informed the author that has recently studied the structure of kangaroo points and that he sees possibilities how to define an invariant which does not increase. Again, one has to wait until written material becomes available.
\med
\med


{\bf References}\med

{\parindent 1.2cm

\litem {[Ab1]} Abhyankar, S.: Local uniformization of algebraic surfaces over ground fields of characteristic $p$. Ann. of Math. 63 (1956), 491-526. 

\litem {[Ab2]} Abhyankar, S.:  Good points of a hypersurface. Adv. Math. 68
(1988), 87-256.

\litem {[BM]} Bierstone, E., Milman, P.: Canonical desingularization in characteristic zero
by blowing up the maximum strata of a local invariant. Invent.÷ Math.÷ 128 (1997), 207-302.

\litem {[BV]} Bravo, A., Villamayor, O.: Hypersurface singularities in positive characteristic and stratification of the singular locus. Preprint 2008.

\litem {[CP]} Cossart, V., Piltant, O.: Resolution of singularities of threefolds in positive characteristic I and II. hal.archives-ouvertes.fr/hal-00139124/fr, hal-00139445/fr.

\litem {[Cu]}  Cutkosky, Resolution of singularities for 3-folds in positive characteristic. arXiv: math.AG/0606530.

\litem{[EiH]} Eisenbud, D., Harris, J.: The geometry of schemes. Graduate Texts in Mathematics. Springer 2000.

\litem {[EH]} Encinas, S., Hauser, H.: Strong resolution of singularities in
characteristic zero. Comment. Math. Helv. { 77} (2002),  421-445.

\litem {[EV1]} Encinas, S., Villamayor, O.: A course on constructive desingularization and equivariance. In:   Resolution of singularities, Obergurgl 1997, Progress in Mathematics 181 (2000).

\litem{[EV2]} Encinas, S., Villamayor, O.: Good points and constructive resolution of singularities. Acta Math. { 181} (1998), 109-158.

\litem {[EV3]} Encinas, S., Villamayor, O.: Rees algebras and resolution of singularities.  Rev. Mat. Iberoamericana. Proceedings XVI-Coloquio Latinoamericano de \'Algebra 
2006.

\litem{[Ha1]} Hauser, H.: Why Hironaka's proof of resolution fails in positive characteristic. Manuscript 2003, available at www.hh.hauser.cc.

\litem{[Ha2]} Hauser, H.: The Hironaka Theorem on resolution of singularities. (Or: A proof that we always wanted to understand.) Bull. Amer. Math. Soc.  40 (2003), 323-403.

\litem{[Ha3]} Hauser, H.:  Excellent surfaces over a field and their taut resolution.  In: Resolution of Singularities, Progress in Math.
181, Birkh\"auser 2000.

\litem{[Ha4]} Hauser, H.: Three power series techniques. Proc.÷ London Math.÷ Soc. 88 (2004), 1-24.

\litem{[Ha5]} Hauser, H.: Seven short stories on blowups and resolution.  In: Proceedings of G\"okova Geometry-Topology Conference 2005 (ed.÷ S.÷ Akbulut et al.),  1-48, International Press, 2006. 

\litem{[Hi1]} Hironaka, H.: Program for resolution of singularities in  characteristics $p > 0$. Notes from lectures at the Clay Mathematics Institute, September 2008.

\litem{[Hi2]} Hironaka, H.: A program for resolution of singularities, in all characteristics $p > 0$ and in all dimensions. Lecture notes ICTP Trieste, June 2006.

\litem{[Hi3]} Hironaka, H.: Theory of infinitely near singular points. J. Korean Math. Soc. 40 (2003), 901-920.

\litem{[Hi4]} Hironaka, H.: Desingularization of excellent surfaces. Notes by B. Bennett at the Conference on Algebraic Geometry, Bowdoin 1967. Reprinted in: Cossart, V., Giraud, J., Orbanz, U.: Resolution of surface singularities. Lecture Notes in Math. 1101, Springer 1984. 

\litem{[Hi5]} Hironaka, H.:  Resolution of singularities of an algebraic variety over a field of characteristic zero.  Ann. of Math. 79 (1964), 109-326.

\litem{[HWZ]} Hauser, H., Wagner, D., Zeillinger, D.: Surface resolution in positive characteristic (using characteristic $0$ invariants). In preparation.

\litem{[MK]} Matsuki, K., Kawanoue, H.: Toward resolution of singularities over a field of positive characteristic (The
Kawanoue program) Part II. Basic invariants associated to the idealistic filtration and their properties.

\litem{[Ka]} Kawanoue, H.: Toward resolution of singularities over a field of positive characteristic. arXiv:math.AG/0607009.

\litem {[Ko]} Koll\'ar, J.: Lectures on resolution of singularities.  Annals Math. Studies, vol. 166. Princeton 2007.

\litem {[Mo]} Moh, T.-T.: On a stability theorem for local uniformization in
characteristic $p$. Publ. Res. Inst. Math. Sci.   23 (1987), 965-973.

\litem {[Mu]} Mulay, S.: Equimultiplicity and hyperplanarity. Proc.÷  Amer.÷  Math.÷ Soc.÷  { 87} (1983), 407-413.

\litem {[Na1]} Narasimhan, R.: Monomial equimultiple curves in positve characteristic.
Proc. Amer. Math. Soc. {89} (1983), 402-413.

\litem {[Na2]} Narasimhan, R.: Hyperplanarity of the equimultiple locus. Proc. Amer.
Math. Soc. { 87} (1983), 403-406.

\litem {[Se]} Seidenberg, A.: Reduction of singularities of the differential equation $A\,dy=B\,dx$.  Amer. J. Math.  90  (1968), 248--269.

\litem{[Vi1]} Villamayor, O.: Constructiveness of Hironaka's resolution. Ann. Sci. \'Ecole Norm. Sup. 22 (1989), 1-32.

\litem{[Vi2]} Villamayor, O.: Hypersurface singularities in positive characteristic. Adv. Math. 213 (2007), 687-733.

\litem{[Vi3]} Villamayor, O.:  Elimination with aplications to singularities in positive characteristic. Publ. Res. Inst. Math. Sci.  44  (2008), 661--697.

\litem {[W\l]} W\l odarczyk, J.: Simple Hironaka resolution in characteristic zero.  J. Amer. Math. Soc.  18  (2005),  779--822.

\litem {[Za]} Zariski, O.: Reduction of the singularities of algebraic three dimensional varieties. Ann. of Math. 45 (1944), 472Ð542.

\litem {[Ze]} Zeillinger, D.: Polyederspiele und Aufl\"osen von Singularit\"aten. PhD Thesis, Universit\"at Innsbruck, 2005.

}

\vs 1cm

Fakult\"at f\"ur Mathematik\par 
Universit\"at Wien, Austria\par
herwig.hauser@univie.ac.at\par

\vfill\eject\end{document}